\documentclass[12pt]{amsart}
\theoremstyle{plain}
\newtheorem{Thm}{Theorem}

\newtheorem{Prop}[Thm]{Proposition}

\begin{document}

\title[ partial convexity]
{Partial convexity to the heat equation}

\author{Li MA }

\address{Li Ma:\ Department of mathematical sciences,
Tsinghua University, Beijing 100084, China.\ \ \
lma@math.tsinghua.edu.cn}

\dedicatory{}
\date{March. 20th, 2009}

\thanks{$^*$ The research is partially supported by the National Natural Science
Foundation of China 10631020 and SRFDP 20060003002. He would like
to thank Prof.Jiaping Wang for helpful discussion. }

\begin{abstract}
In this paper, we study the partial convexity of smooth solutions
to the heat equation on a compact or complete non-compact
Riemannian manifold $M$ or on Kaehler-Ricci flow. We show that
under a natural assumption, a new partial convexity property for
smooth solutions to the heat equation is preserved.
\end{abstract}

 \maketitle

\section{Introduction}
It is always interesting for people to find invariant sets for
evolution equations. Invariant sets of geometric evolution
equations can be defined by related monotonicity quantities. In
the recent study of Ricci flow, G.Perelman \cite{P02} finds new
monotonicity formulae. In this paper, we study invariant sets of
the heat flow by the maximum principle trick. More precisely, we
are concerned with partial convexity of solutions to the heat
equation (associated to Ricci flow or Ricci-Kahler flow)
 by using the maximum principle method.
 Consider a smooth solution to the heat equation
\begin{equation}
u_t=\Delta u,\;\;\;on \;\;M\times[0,T) \label{ht}
\end{equation}
on the compact or complete Riemannian manifold $(M^n,g)$ of
dimension $n$. Here $\Delta$ is the Laplacian operator of the
metric $g$ (with the sign such that $\Delta u=u^{''}$ on the real
line $R$ ). By definition, for some
 positive integer
 $1\leq k\leq n$, we say that
 $u$ is \emph{partially k-convex} on $M\times [0,T)$ if  we have
 that the positivity of the
 function
 $$\sigma_j(u)>0, \quad{j=1,...,k,}$$
 where $\sigma_j(u)$ is the j-th elementary symmetric
 polynomial of the hessian matrix $D^2u$ on $(M,g)$), of the solution $u$
 is preserved along the heat equation (\ref{ht}).
 We also study the corresponding problem for heat equation associated to the Kaehler-Ricci flow.

 Let's just mention some earlier related research.
 In the
 previous study, R.Hamilton \cite{H93} has extended the famous Li-Yau \cite{LY86}
 gradient estimate
 of the heat equation and has obtained the matrix Harnack inequality
 for the heat equation on Riemannian manifolds. Cao and Ni \cite{CN03} have studied
 the  matrix Harnack inequality
 for the heat equation along the Kahler-Ricci flow. Some years ago,
 Brascamp and E.Lieb etc \cite{Lieb}
 have studied the log convexity of the solution to the heat equation on convex
 domains, and soon after, S.T.Yau has found the continuity
 argument for convexity of solution (see the appendix in \cite{Y},
 also \cite{Kore82},
 \cite{Y2} for this point and we shall use this kind of argument).
  L.Caffarelli and A.Friedman \cite{CaF85}, P.L.Lions\cite{Lions04}, etc, have studied
 the convexity of more general elliptic and parabolic equations.
 We refer to the works of Korevaar \cite{Kore82} and Kawohl \cite{Ko85} for more results
 and references.

To state our results, we make some conventions. Let $A=(u_{ij})$ be
the Hessian matrix of the smooth function $u$. Here, we have denoted
by $\nabla_i f=f_i$, in local coordinates $(x^i)$, the i-th
covariant derivative of a smooth function $f$ on $M$, and $f_{ij}$
the corresponding second covariant derivatives.

Take a fixed point $x\in M$ and normal coordinates $(x^i)$ at
$x\in M$. We now consider some elementary algebra on the tangent
space $T_xM$. Note that at this point $x$ we have $(g_{ij})=I$,
which is the identity. For $k=1,2,...,n$, let $\sigma_k(A)$ be the
k-th elementary symmetric polynomial of $A$. For example, $$
\sigma_2(A)=\sum_{i<j}\lambda_i\lambda_j
$$
for the diagonal matrix $A=(\lambda_1\bigoplus...
\bigoplus\lambda_n)$. The $k$-th Newton transformation associated
to $A$ is
$$
T_k(A)=\sigma_k(A)I-\sigma_{k-1}(A)A+...+(-1)^kA^k.
$$
In particular, we have $$ T_1(A)=\sigma_1(A)I-A.
$$
Note that $\sigma_1(A)=\Delta u$.

 Let $A(s)$ be a smooth one-parameter family of symmetric matrices on $T_xM$.
Then we have
$$
\frac{d}{ds}\sigma_k(A(s))=trace(T_{k-1}(A(s)\circ
\frac{d}{ds}A(s)).
$$
Hence, we have
\begin{eqnarray*}
\frac{d}{ds}\sigma_2(A(s))&=&trace(T_{1}(A(s))\circ
\frac{d}{ds}A(s))\\
&=&\sigma_1(A)\frac{d}{ds}\sigma_1(A)-trace(A(s)\circ
\frac{d}{ds}A(s))
\end{eqnarray*}

 One may see \cite{Rei73} for more relations. Sometimes,
 we also denote by $\sigma_k(u)=\sigma_k(A)$. Clearly, one can define
 similar concepts on Kaehler manifolds.

 For the heat equation associated to the Kaelher-Ricci flow,
 we have the following result
\begin{Thm}\label{thm1} Let $(M,g(t))$ be a compact or complete non-compact Kaehler manifold
of dimension $n$, where $(g(t)), 0\leq t<T,$ is a Kaehler-Ricci
flow with bounded curvature. Assume that the holomorphic
bi-sectional curvature of each $g(t)$ is non-negative such that
\begin{equation}\label{non1}
-u_{\beta\bar{\alpha}}R_{\alpha\bar{\beta}\gamma\bar{\delta}}u_{\bar{\gamma}\delta}+
u_{\beta\bar{\alpha}}R_{\alpha\bar{s}}u_{s\bar{\beta}}\geq 0,
\end{equation}
for any hermitian matrix $(u_{\bar{\gamma}\delta})$. Let
$\Delta_{g(t)}$ be the Laplacian of $g(t)$, $0\leq t<T$. Let $u$
be a smooth solution to
$$
u_t=\Delta_{g(t)} u
$$ on $M^n\times[0,T)$ with \emph{nice decay at infinity} when $M$ noncompact.
Let $A=(u_{\alpha\bar{\beta}})$ be the hessian of $u$ with
$T_{\alpha}$ and $\sigma_{\alpha}$, $\alpha=1,2$, defined as
above. Then we have
$$
\sigma_1(A)\geq (>)0, \;\;\;for\;\;t>0
$$
provided $\sigma_1(A)\geq (>)0$ at $t=0$. Furthermore, we have
partial 2-convexity of the solution $u$; that is, the positivity
of the function
 $\sigma_2(u)$
 is also preserved provided it is positive at $t=0$.
\end{Thm}

By definition, a smooth function $f\in C^2(M)$ has a \emph{nice
decay at infinity} if $|f|(x)+|\nabla f|(x)+|D^2f(x)|\to 0$ as the
distance $d(x,0)\to \infty$ for some fixed point $o\in M$. We recall
here that $g(t)$ is a Kaehler-Ricci flow on the manifold $M$ if in
local complex coordinates $(z^{\alpha})$, we have
$$
\partial_tg_{\alpha\bar{\beta}}=-R_{\alpha\bar{\beta}},
\;\;on\;\;M\times (0,T).
$$

We remark that in the proof of Theorem \ref{thm1} in section
\ref{sect2}, we only use condition \ref{non1} for any hermitian
symmetric matrix $(u_{ij})$ in place of the non-negative
holomorphic bi-sectional curvature. In particular, the assumption
(\ref{non2}) is automatically true on standard n dimensional
complex projective space $CP^n$. We believe that one may apply
Theorem \ref{thm1} to the
 study of Kaehler-Ricci flow.

  To explain the idea in the proof of Theorem \ref{thm1},
let's assume that we are studying periodic or nice decay solutions
to the heat equation on $R^n$. In this case, we have
$$
A_t=\Delta A.
$$
Then we have
\begin{equation}
\sigma_1(A)_t=\Delta \sigma_1(A).\label{sigma1}
\end{equation}
Hence, assuming $\sigma_1(A)>0$ at $t=0$, by the standard maximum
principle, the positivity $\sigma_1(u)$ is preserved along the
heat equation.
 Here we have assumed that $u(.,t)$ has a nice decay at infinity for
each $t\in [0,T)$. In fact, this has be proved in \cite{DM}. We
may make similar computation at least formally for higher order
symmetric elementary functions. Since,
$$
\partial_t\sigma_k(A)=trace(T_{k-1}(A)\circ \partial_tA),
$$
see (2.3) in \cite{GV}, and
$$
\Delta \sigma_k(A)=trace(\nabla_iT_{k-1}(A)\circ
\nabla_iA)+trace(T_{k-1}(A)\circ \Delta A),
$$
we have
\begin{eqnarray*}
(\partial_t-\Delta)\sigma_k(A)&=&-trace(\nabla_iT_{k-1}(A)\circ
\nabla_iA)+trace(T_{k-1}(A)\circ (\partial_t-\Delta) A)\\
&=&
-trace(\nabla_iT_{k-1}(A)\circ \nabla_iA).
\end{eqnarray*}
Then by the maximum principle, we have the following

\begin{Prop}\label{thm2} Let $u$ be a smooth solution to (\ref{ht}) on $R^n\times[0,T)$
with nice decay at infinity. Let $A=(u_{ij})$ be the hessian of $u$
with $T_k$ and $\sigma_k$ defined above. Assume that
\begin{equation} trace(\nabla_iT_{k-1}(A)\circ
\nabla_iA)\leq 0,\;\;on\;\;M\times[0,T).\label{ass}
\end{equation}
Then we have
$$
\sigma_k(A)\geq (>)0, \;\;\;for\;\;t>0
$$
provided $\sigma_k(A)\geq (>)0$ at $t=0$.
\end{Prop}

 We remark that the assumption (\ref{ass}) may be difficult to verify in applications.
 However, we can find some partial convexity condition for the heat equation.
 In fact, we have
\begin{eqnarray*}
(\partial_t-\Delta)\sigma_2(A)&=&-trace(\nabla_iT_{1}(A)\circ
\nabla_iA)\\
&=&-trace(\nabla_i\sigma_{1}(A)\nabla_iA)+trace(\nabla_i
A\circ \nabla_iA)\\
&=&-|\nabla\sigma_{1}(A)|^2+trace(\nabla_i A\circ \nabla_iA)\\
&\geq& -|\nabla\sigma_{1}(A)|^2.
\end{eqnarray*}
By (\ref{sigma1}), we have
$$
(\partial_t-\Delta)\sigma_1(A)^2/2=-|\nabla\sigma_1(A)|^2.
$$
Set $F=\sigma_1(A)^2/2$ and $H=\sigma_2(A)/F$.
 Then, using the fact that $\sigma_1^2\geq 2\sigma_2$
 and the formula
 $$
L(\frac{u}{v})=\frac{Lu}{v}-\frac{uLv}{v^2}-2\frac{u}{v^3}|\nabla
v|^2+\frac{2}{v^2}(\nabla u,\nabla v)
 $$
 where $L=\partial_t-\Delta$,
 we obtain that
 $$
(\partial_t-\Delta)H\geq 2<\nabla H,\nabla F>/F,
$$
 and we get by the maximum principle that
$$
\sigma_2(A)/F\geq (>)0, \;\;\;for\;\;t>0
$$
provided $\sigma_2(A)/F\geq (>)0$ at $t=0$. Hence, assuming
$\sigma_1(A)>0$, the positivity property of $\sigma_2(A)$ is
preserved along the heat equation. Hence, we have

\begin{Prop}\label{thm0} Let $u$ be a smooth solution to (\ref{ht}) on $R^n\times[0,T)$
with nice decay at infinity. Let $A=(u_{ij})$ be the hessian of
$u$ with $T_k$ and $\sigma_k$ defined above. Assume that
$\sigma(A)>0$ and $\sigma_1(A)>0$ at $t=0$. Then we have
$$
\sigma_2(A)\geq (>)0, \;\;\;for\;\;t>0.
$$
\end{Prop}

For the heat equation in the Riemannian case when $g$ being a fixed
metric, we have the following result whose assumption is similar to
Corollary 4.4 in \cite{H93}.

\begin{Thm}\label{thm3} Let $(M,g)$ be a compact or complete noncompact Riemannian manifold
of dimension $n$ with non-negative sectional curvature and
parallel Ricci curvature tensor. Assume further that
\begin{equation}\label{non2} -2u_{ij}R_{ikjl}u_{kl}+2u_{ij}R_{jl}u_{il}\geq
0, \; on \; M,
\end{equation}
for any symmetric matrix $(u_{ij}$. Let $\Delta$ be the Laplacian
of $g$. Let $u$ be a smooth solution to
$$
u_t=\Delta u
$$ on $M^n\times[0,T)$ with nice decay at infinity when $M$ is complete and noncompact.
Let $A=(u_{ij})$
be the hessian of $u$ with $T_k$ and $\sigma_k$ defined above. Then
we have
$$
\sigma_1(A)\geq (>)0, \;\;\;for\;\;t>0
$$
provided $\sigma_1(A)\geq (>)0$ at $t=0$. Furthermore, we have
partial 2-convexity of the solution $u$; that is, the positivity of
the function
 $$\sigma_2(u)$$
 is preserved provided it is positive at $t=0$.
\end{Thm}
We remark that in the proof of Theorem \ref{thm3}, we only use
condition \ref{non2} for any symmetric matrix $(u_{ij})$, not the
non-negative section curvature. In particular, the assumption
(\ref{non2}) is automatically true on standard n-sphere $S^n$.

\section{proof of Theorem \ref{thm1}}\label{sect2}

Let $(M^m,g(t))$ be a compact Kaehler manifold, where $g(t)$ is a
Kaehler-Ricci flow in the sense that
$$
\partial_tg_{\alpha\bar{\beta}}=-R_{\alpha\bar{\beta}}.
$$

Let $\Delta=\Delta_{g(t)}$ be the Laplacian of the metric $g(t)$.
Assume that $u\in C^2(M\times[0,T))$ satisfies the heat equation
$$
(\partial_t-\Delta)u=0,\;\;\;on\;\;M_T
$$
where $M_T=M\times[0,T)$. Doing the computation as in Lemma 2.1 in
(\cite{NT03}), we have
$$
(\partial_t-\Delta)u_{\alpha\bar{\beta}}=R_{\alpha\bar{\beta}\gamma\bar{\delta}}
u_{\bar{\gamma}{\delta}}
-\frac{1}{2}(R_{\alpha\bar{s}}u_{s\bar{\beta}}+u_{\alpha\bar{s}}R_{s\bar{\beta}})
$$
Let $A=(u_{\alpha\bar{\beta}})$. Then we have
$$ (\partial_t-\Delta)trace(A)=0,
$$
which implies that $\sigma_1(0)>0 (\geq 0)$ is preserved on the heat
equation,
 and then $$
trace(A\circ(\partial_t-\Delta)A)=
u_{\beta\bar{\alpha}}R_{\alpha\bar{\beta}\gamma\bar{\delta}}u_{\bar{\gamma}\delta}-
u_{\beta\bar{\alpha}}R_{\alpha\bar{s}}u_{s\bar{\beta}}.
$$
Hence, we have
\begin{eqnarray*}
&&(\partial_t-\Delta)\sigma_2(A)\\&=&-trace(\nabla_{\alpha}T_{1}(A)\circ
\nabla_{\bar{\alpha}}A)+trace(T_1(A)\circ (\partial_t-\Delta) A)\\
&=& -|\nabla \sigma_1(A)|^2+trace(\nabla_{\alpha}(A)\circ
\nabla_{\bar{\alpha}}A)+\sigma_1(A)trace((\partial_t-\Delta) A)\\
&-&trace(A\circ(\partial_t-\Delta) A)\\
&=&-|\nabla \sigma_1(A)|^2+trace(\nabla_{\alpha}(A)\circ
\nabla_{\bar{\alpha}}A)-
u_{\beta\bar{\alpha}}R_{\alpha\bar{\beta}\gamma\bar{\delta}}u_{\bar{\gamma}\delta}+
u_{\beta\bar{\alpha}}R_{\alpha\bar{s}}u_{s\bar{\beta}}.
\end{eqnarray*}
Note that by our assumption (\ref{non1}), $$ -
u_{\beta\bar{\alpha}}R_{\alpha\bar{\beta}\gamma\bar{\delta}}u_{\bar{\gamma}\delta}+
u_{\beta\bar{\alpha}}R_{\alpha\bar{s}}u_{s\bar{\beta}}\geq 0. $$
Then, using the same trick as in what we did in Proposition
\ref{thm0}, we get the partial convexity for $\sigma_2$ by the
maximum principle.

\section{proof of Theorem \ref{thm3}}\label{sect3}

Recall the second contracted Bianchi identity that
$$
\nabla_kR_{ikkl}+\nabla_jR_{il}-\nabla_lR_{ij}=0.
$$
Using the Ricci formula and the identity above, we can compute that
\begin{eqnarray}
&&\nabla_i\nabla_j(\Delta
u) \label{ric}\\
&=&\Delta\nabla_i\nabla_ju+2R_{kijl}\nabla_k\nabla_lu-R_{il}\nabla_j\nabla_lu-
R_{jl}\nabla_i\nabla_lu \nonumber\\
&&-(\nabla_iR_{jl}+\nabla_jR_{il}-\nabla_lR_{ij})\nabla_lu.
\nonumber
\end{eqnarray}
In fact, we have
\begin{eqnarray*}
&&\nabla_i\nabla_j(\Delta u)=u_{kkji}\\
&=& \nabla_i(u_{jkk}-R_{jl}u_l)\\
&=&u_{jkki}-\nabla_iR_{jl}u_l-R_{jl}u_{li}\\
&=&u_{jkik}-R_{ikjm}u_{km}-R_{im}u_{mj}-\nabla_iR_{jl}u_l-R_{jl}u_{li}\\
&=&\nabla_k(u_{jik}-R_{ikjm}u_m)-R_{ikjm}u_{km}
-R_{im}u_{mj}-\nabla_iR_{jl}u_l-R_{jl}u_{li}\\
&=&u_{ijkk}-\nabla_kR_{ikjm}u_m-R_{ikjm}u_{mk}\\
&- & R_{ikjm}u_{km}
-R_{im}u_{mj}-\nabla_iR_{jl}u_l-R_{jl}u_{li}.
\end{eqnarray*}
So, we get (\ref{ric}) by using the second Bianchi identity above.

 By our assumption that
$$
\nabla_jR_{il}=0,
$$
we have
$$
\nabla_i\nabla_j(\Delta
u)=\Delta\nabla_i\nabla_ju+2R_{kijl}\nabla_k\nabla_lu-R_{il}\nabla_j\nabla_lu-
R_{jl}\nabla_i\nabla_lu.
$$
Hence we have
$$
\partial_tu_{ij}=\Delta u_{ij}+2R_{kijl}u_{kl}-R_{il}u_{jl}-
R_{jl}u_{il}.
$$
Set $A=(u_{ij})$. Again we have that $$ trace((\partial_t-\Delta)
A)=0,
$$
which also gives us that $\sigma_1(u)>0$ is preserved on the heat
equation. We now have that
\begin{eqnarray*}
&&(\partial_t-\Delta)\sigma_2(A)\\&=&-trace(\nabla_{i}T_{1}(A)\circ
\nabla_{i}A)+trace(T_1(A)\circ (\partial_t-\Delta) A)\\
&=& -|\nabla \sigma_1(A)|^2+trace(\nabla_i(A)\circ
\nabla_{i}A)+\sigma_1(A)trace((\partial_t-\Delta) A)\\
&-&trace(A\circ(\partial_t-\Delta) A)\\
&=&-|\nabla \sigma_1(A)|^2+trace(\nabla_{i}(A)\circ \nabla_{i}A)-
2u_{ij}R_{ikjl}u_{kl}+2u_{ij}R_{jl}u_{il}.
\end{eqnarray*}
Note that by our assumption (\ref{non2})), we always have that $$
-2u_{ij}R_{ikjl}u_{kl}+2u_{ij}R_{jl}u_{il}\geq 0.
$$
Hence, we have $$ (\partial_t-\Delta)\sigma_2(A)\geq -|\nabla
\sigma_1(A)|^2.
$$
 Then, using again the same trick as in what we did in Proposition \ref{thm0}, we get the desired partial
convexity by the maximum principle.

\end{document}